\documentclass{amsart}

\newtheorem{lemma}{Lemma}[section]
\newtheorem{theorem}[lemma]{Theorem}

\newcommand{\intbar}{\mathop{\int\makebox(-13.5,0){\rule[4pt]{.7em}{0.3pt}}
\kern-6pt}\nolimits} 

\makeatletter
\def\proof{\@ifnextchar[\opargproof{\opargproof[\bf Proof. ]}}
\def\opargproof[#1]{\par\noindent {\bf #1}}

\makeatother
\def\ep{\varepsilon}
\def\al{\alpha}
\def\Om{\Omega}
\def\om{\omega}
\def\lf{\left}
\def\rg{\right}
\def\p{\partial}
\def\ds{\displaystyle}

\newcommand{\R}{\mathbb{R}}
\newcommand{\N}{\mathbb{N}}

\begin{document}

\title{Partial regularity for harmonic maps, and related problems}
\author{Tristan Rivi\`ere}
\address[T.~Rivi\`ere]{Mathematik\\ETH-Z\"urich\\CH-8092 Z\"urich}
\email{riviere@math.ethz.ch}
\author{Michael Struwe}
\address[M.~Struwe]{Mathematik\\ETH-Z\"urich\\CH-8092 Z\"urich}
\email{struwe@math.ethz.ch}

\date{\today}

\begin{abstract}
Via gauge theory, we give a new proof of partial regularity for harmonic maps 
in dimensions $m \ge 3$ into arbitrary targets. This proof avoids the use of adapted frames
and permits to consider targets of "minimal" $C^2$ regularity. The proof we
present  moreover extends to a large class of elliptic systems
of quadratic growth.
\end{abstract}

\maketitle

\section{Introduction}
In \cite{Riviere}, the first author presented a new approach to the regularity 
result of H\'elein \cite{Hel} for weakly harmonic maps in dimension $m=2$ where he 
succeeded in writing the harmonic map system in the form of a conservation law
whose constituents satisfied elliptic equations with a Jacobian structure to 
which Wente's \cite{Wente} regularity results could be applied. 

Consider for instance a harmonic map $u=(u^1,\dots,u^n)\in H^1(B;\R^n)$
from a ball $B^m = B \subset \R^m$ to a hypersurface $N\subset \R^n$ with normal $\nu$. 
In this case the harmonic map equation may be written in the form 
\begin{equation}\label{1.1}
  -\Delta u^i= w^i \nabla w^j \cdot \nabla u^j = 
  (w^i \nabla w^j -w^j \nabla w^i) \cdot \nabla u^j,\, 1\le i \le n,
\end{equation}
where $w = \nu \circ u$. The key idea then is to identify the anti-symmetry of the $1$-form 
\begin{equation}\label{1.2}
  \Omega^{ij} = (w^i d w^j -w^j d w^i),\, 1\le i,j \le n,
\end{equation}
as the essential structure of equation \eqref{1.1}. 

Interpreting $\Omega \in L^2(B; so(n)\otimes \wedge^1{\R}^n)$ as a connection in the
$SO(n)$-bundle $u^* TN$ and following Uhlenbeck's approach the existence of Coulomb gauges
\cite{Uhl}, one succeeds in finding $P \in H^1(B;SO(n))$ and $\xi \in H^1(B)$ 
such that 
\begin{equation}\label{1.3}
  P^{-1}dP +P^{-1}\Omega P = *d\xi,
\end{equation}
where $*$ is the Hodge dual. If $m = 2$, further algebraic manipulations then yield 
the existence of $A, B \in H^1(B)$ with 
\begin{equation}\label{1.4}
     ||dist(A,SO(n))||_{L^{\infty}} \le C||\Omega||_{L^2}
\end{equation}
such that \eqref{1.1} may be written as 
\begin{equation}\label{1.5}
   div(A\nabla u+B\nabla^\perp u)=0,
\end{equation}
where $\nabla^\perp = *d$. 
By Hodge decomposition one then obtains $E$ and $D$ in $W^{1,2}(D^2)$ such that
\begin{equation}\label{1.6}
   A\,\nabla u=\nabla D + \nabla ^\perp E.
\end{equation}
From \eqref{1.5} we see that $D$ and $E$ satisfy the equations
\begin{equation}\label{1.7}
  \begin{split}
    & -\Delta D = div(A\,\nabla u)=-\nabla B\cdot\nabla^\perp u,\\
    & -\Delta E = curl(A\,\nabla u)=\nabla^\perp A\cdot\nabla u,
  \end{split}
\end{equation}
which exhibit the desired Jacobian structure. The results in \cite{CLMS} then
imply that $D, E \in W_{loc}^{2,1}(B)$. Provided that
we restrict our attention to a domain where $||\Omega||_{L^2}$ is sufficiently small,
from \eqref{1.4} we conclude that 
$\nabla u=A^{-1}(\nabla D +\nabla^\perp E)\in W_{loc}^{1,1}(B)$ 
and $u \in W_{loc}^{2,1}(B) \hookrightarrow C^0(B)$, which implies full regularity. 

In dimensions $m \ge 3$, the harmonic map equation is super-critical in the Sobolev space
$H^1(B;{\R}^n)$ and no regularity result, not even a partial one, can be expected. 
In fact, in \cite{Ri1} the first author constructed examples of weak solutions to 
\eqref{1.1} in $H^1(B;S^2)$ for $m \ge 3$ which are {\it nowhere}
continuous.  

Under the further assumption that the solution $u$ lies in the homogeneous Morrey space 
$L^{2,m-2}_1$, which sometimes also is denoted as $M^{1,2}_2$, with 
\begin{equation}
\label{00}
\|u\|^2_{L^{2,m-2}_1}=\sup_{x\in B\, ,\, r>0}\left(\frac{1}{r^{m-2}}\int_{B_r(x)\cap B}|\nabla u|^2\right)<+\infty,
\quad,
\end{equation}
the harmonic map equation (\ref{1.1}) becomes critical. More generally, this is true for
any elliptic system with a nonlinearity growing quadratically in the gradient (see \cite{Gia}).
Assumption \eqref{00} is natural in the context of harmonic maps; in fact, it is a 
direct consequence of the geometric stationarity assumption described in Section 2. 
Observe that in dimension $m=2$ assumption \eqref{00} corresponds exactly to the assumption
of finite energy and it therefore appears as the natural extension of the finite energy
condition to higher dimensions. 

Strengthening the assumption that $u\in H^1(B;\R^n)$ by assuming \eqref{00}, one might then 
hope to be able to extend the approach described above to the case $m\ge 3$. 
However, in order to achieve \eqref{1.4}, in dimension $m=2$ one crucially uses that by
Wente's result mentioned above the solution $\psi \in H^1(B^2;\wedge^2\R^n)$
of the equation 
\begin{equation}
\label{02}
\lf\{
\begin{array}{l}
\ds\Delta\psi=df\wedge dg\quad\quad\mbox{ in }B^m\\[5mm]
\ds \psi=0\quad\quad\mbox{ on }\p B^m
\end{array}
\rg.
\end{equation}
for given $f$ and $g$ in $H^1(B^2)$ belongs to $L^\infty$.
Unfortunately, 
this result does not extend to $m\ge 3$ when we replace the assumption $f,g \in H^1(B^m)$
by the condition that $f$ and $g$ belong to the Morrey space $L^{2,m-2}_1(B^m)$. Indeed, for $m=3$,
letting $f =\frac{x_1}{|x|}$ and choosing $g = \frac{x_2}{|x|}$, we have $f,g \in L^{2,1}_1(B^3)$
and equation \eqref{02} admits a unique solution $\psi \in L^{2,1}_1(B^3)$, but $\psi \notin L^\infty$.
Thus the $L^{\infty}$-bound \eqref{1.4} does not seem to be available in dimension larger 
than $2$ and the approach outlined above seems to fail for this reason.

However, as we presently explain, \eqref{1.1} - \eqref{1.3} in combination with standard techniques 
of elliptic regularity theory already suffice to conclude partial regularity, directly. 
In fact, via the gauge transformation $P$, from \eqref{1.1} we obtain the equation 
\begin{equation}\label{1.8}
  -div(P^{-1}\nabla u) = (P^{-1}\nabla P +P^{-1}\Omega P) \cdot P^{-1}\nabla u 
  = *d\xi \cdot P^{-1}du,
\end{equation}
where the right hand side already has the structure of a Jacobian -- up to the harmless 
(bounded) factor $P^{-1}$. Also observe that $\nabla u$ may be recovered from the 
term $P^{-1}\nabla u$ without any difficulty. 

More generally, partial regularity results can be obtained for a large class of elliptic 
systems with quadratic growth that can be cast in the form 
\begin{equation}
\label{04}
-\Delta u=\Om\cdot\nabla u\quad\quad\mbox{ in }B
\end{equation}
already considered in \cite{Riviere}. (In coordinates, equation \eqref{04} simply reads
$-\Delta u^i=\Om^{ij}\cdot\nabla u^j$.)

\begin{theorem}
\label{th-I.1}
For every $m\in\N$ there exists $\ep(m)>0$ such that for every 
$\Om \in L^2(B^m,so(n)\otimes \wedge^1{\R}^m)$ and for every weak solution 
$u \in H^1(B^m,{\R}^n)$ of equation \eqref{04}, satisfying 
the Morrey growth assumption
\begin{equation}
\label{05}
  \sup_{x\in B,\, r>0}\left(\frac{1}{r^{m-2}}
  \int_{B_r(x)\cap B}(|\nabla u|^2+|\Om|^2)\, dx \right)<\ep(m)\quad,
\end{equation}
we have that $u$ is locally H\"older continuous in $B$ with exponent $0<\al = \al(m) <1$.
\end{theorem}

The previous result is optimal, as shown by the standard example of the weakly harmonic map
$u \colon B^3 \to S^2 \hookrightarrow {\R}^3$ with $u(x)=x/|x|$.
We have $u \in H^1(B^3,{\R}^3)$ and, letting 
$\Om=(\Om^{ij}):=(u^id u^j-u^j du^i)\in L^2(B^m,so(n)\otimes \wedge^1{\R}^m)$,
we see that $u$ weakly satisfies the equation (\ref{04}) and the condition
\begin{equation}
\label{06}
\sup_{x\in B\, ,\, r>0}\left(\frac{1}{r^{m-2}}\int_{B_r(x)\cap B}(|\nabla u|^2+|\Om|^2)\, dx \right)<+\infty\quad.
\end{equation}
The map $u$, however, is not continuous at the origin.

\section{Stationary harmonic maps}
For a smooth, compact, oriented $k$-dimensional submanifold $N \subset \R^n$ 
and a ball $B \subset \R^m$ let
\begin{equation}
  H^1(B;N) = \{u \in H^1(B;\R^n); \, u(x) \in N \hbox{ for almost every } 
  x \in B\}. 
\end{equation}
Recall that a map $u\in H^1(B;N)$ is {\em stationary} if $u$ is critical for the energy 
\begin{equation*}
    E(u)=\int_{B}|\nabla u|^2\, dx
\end{equation*}
both with respect to variations of the map $u$ and with respect to variations in the 
domain. 

It follows that $u$ is weakly harmonic; that is, $u$ satisfies the equation
\begin{equation}\label{2.1}
  -\Delta u= A(u)(\nabla u,\nabla u)
  = \sum_{l=1}^{n-k}\sum_{\alpha=1}^{m}\nu_l \langle d \nu_l \partial_{\alpha}u,
     \partial_{\alpha} u\rangle)
  = \sum_{l=1}^{n-k}w_l \langle \nabla w_l,\nabla u\rangle
\end{equation} 
in the sense of distributions, where $A$ is the second fundamental form of $N$,
defined via an orthonormal frame field $\nu_l$, $1 \le l \le n-k$ for the normal 
bundle to $N$.
Again we denote as $w_l = \nu_l \circ u$ the corresponding unit normal vector 
field along the map $u$, and we denote as $\langle\cdot,\cdot\rangle$ the 
Euclidean inner product.

Moreover, as a consequence of the stationarity condition with respect to 
variations in the domain we have the monotonicity estimate
\begin{equation}\label{2.2}
     r^{2-m} \int_{B_r(x_0)} |\nabla u|^2\, dx \le R^{2-m} \int_{B_R(x_0)} |\nabla u|^2\, dx
\end{equation}
for all balls $B_R(x_0) \subset B$ and all $r \le R$.

The following result was obtained by Evans \cite {Evans} and Bethuel \cite{Bethuel}. 
Note that their approach in general requires the target manifold $N^k$ to be of class $C^5$; see
\cite{Hel}, Theorem 4.3.1 and Remark 4.3.2.
As a corollary to our main result Theorem~\ref{th-I.1}, however, we now easily obtain the 
following generalization of their result to manifolds of class $C^2$.

\begin{theorem}\label{thm2.1}
Let $N^k \subset \R^n$ be a closed submanifold of class $C^2$. Let $m \ge 3$ and 
suppose $u\in H^1(B^m;N)$ is a stationary harmonic map. 
There exists a constant $\varepsilon_0 > 0$
depending only on $N$ with the following property. Whenever on some ball $B_R(x_0) \subset B$
there holds 
\begin{equation}\label{2.3}
    R^{2-m} \int_{B_R(x_0)} |\nabla u|^2\, dx < \varepsilon_0,
\end{equation}
then $u$ is H\"older continuous (and hence smooth) on $B_{R/2}(x_0)$.
\end{theorem}
\begin{proof}
As in \eqref{1.1}, equation \eqref{2.1} equivalently may be written in the form 
 \begin{equation}\label{3.1}
  -\Delta u^i= \Omega^{ij} \cdot \nabla u^j,
\end{equation}
where $\Omega \in L^2(B; so(n)\times \wedge^1{\R}^n)$ in view of our assumption on $N$,
with components
\begin{equation}\label{3.2}
  \Omega^{ij} =  \Omega^{ij} _{\alpha} dx^{\alpha} = 
  \sum_{l=1}^{n-k}(w_l^i d w_l^j -w_l^j d w_l^i),\, 1\le i,j \le n.
\end{equation}
Note that \eqref{2.2} and \eqref{2.3} imply that $\Omega$ belongs to the Morrey
space $L^{2,m-2}(B)$ with 
\begin{equation}\label{3.3}
   \begin{split}
       ||\Omega||^2_{L^{2,m-2}}
       & = \sup_{x_0\in B} r^{2-m}\int_{B_r(x_0)\cap B} |\Omega|^2\, dx \\
       & \le C\sup_{x_0\in B} r^{2-m}\int_{B_r(x_0)\cap B} |\nabla u|^2\, dx 
       \le C\varepsilon_0.
   \end{split}
\end{equation}
The result now is an immediate consequence of Theorem~\ref{th-I.1}.
\end{proof}

\section{Proof of Theorem~\ref{th-I.1}}
We may assume that condition \eqref{05} is satisfied on $B = B_1(0)$.
As in \eqref{1.3}, we obtain the existence of a suitable 
gauge transformation $\Phi$, 
transforming $\Om$ into Coulomb gauge by applying the following lemma. The bound \eqref{05}
also yields corresponding estimates for $P$ and $\xi$. 

\begin{lemma}
\label{lm-3.1}
There exists $P \in H^1(B;SO(n))$ and $\xi \in H^1(B)$ such that 
\begin{equation}\label{3.4}
  P^{-1}dP +P^{-1}\Omega P = *d\xi \mbox{ on } B, \mbox{ and }\xi=0\mbox{ on }\p B\, .
\end{equation}
Moreover, $dP$ and $d\xi$ belong to $L^{2,m-2}(B)$ with
\begin{equation}\label{3.5}
 \begin{split}
   C||dP||^2_{L^{2,m-2}}  + C||d\xi||^2_{L^{2,m-2}} 
   \le C (||\Omega||^2_{L^{2,m-2}}+||du||^2_{L^{2,m-2}}) \le C\varepsilon(m).
  \end{split}
\end{equation}
\end{lemma}
The proof of this lemma will be given in the next section.

Recall that a function $f \in L^1(B)$ belongs to the space $BMO(B)$ if
\begin{equation*}
       [f]_{BMO}
       = \sup_{x_0 \in B, \, r >0}\left(\intbar_{B_r(x_0)\cap B} |f - \bar{f}_{x_0,r}|\, dx\right) < \infty,
\end{equation*}
where
\begin{equation*}
   \bar{f}_{x_0,r} = \intbar_{B_r(x_0)\cap B} f \, dx
\end{equation*}
denotes the average of $f$ over $B_r(x_0)\cap B$, and so on. 
By Poincar\'e's inequality, moreover,
for $1 \le p \le m$ any function $f \in W^{1,p}(B)$ with $df \in L^{p,m-p}(B)$ 
belongs to $BMO(B)$ and there holds 
\begin{equation*}
\label{zx1}
       [f]^p_{BMO} \le C  ||df||^p_{L^{p,m-p}} 
       = \sup_{x_0, \, r >0}\left( r^{p-m}\int_{B_r(x_0)\cap B} |df|^p\, dx\right).
\end{equation*}

Applying the gauge transformation $P^{-1}$ to $\nabla u$ and observing the identity 
$dP^{-1}=- P^{-1}dP P^{-1}$, from \eqref{04} we obtain the equation \eqref{1.8},
that is 
\begin{equation}\label{3.6}
  -div(P^{-1}\nabla u) = (P^{-1}\nabla P +P^{-1}\Omega P) \cdot P^{-1}\nabla u 
  = *d\xi \cdot P^{-1}du.
\end{equation}

Fix a smooth cut-off function $\tau \in C^{\infty}_0(B)$ such that $0 \le \tau \le 1$,
$\tau = 1$ on $B_{1/2}(0)$. Multiplying \eqref{3.6} by $\tau$, we obtain the equation
\begin{equation}\label{3.6a}
  -div(P^{-1}\nabla (u \tau)) 
  = *d\xi \cdot P^{-1}d(u \tau) - e,
\end{equation}
with ``error'' term
\begin{equation}\label{3.6b}
  e = div(P^{-1} u \nabla \tau)) + P^{-1}\nabla u \cdot \nabla \tau +
  *d\xi \cdot P^{-1}u d \tau.
\end{equation}
Since $u\in L^{2,m-2}_1(B)$ and since $u$ is in $L^2(B)$, because of (\ref{zx1}), 
$u$ is in $L^p(B)$ for every $p< \infty$. Therefore, a direct application
of H\"older inequality tells us that $u$ is in $L^{p,m-\delta}$ for every $p<\infty$ 
and for every $\delta>0$. Using this last observation and the fact that
$d\xi$ is in $L^{2,m-2}$, we conclude that 
\begin{equation}
\label{zx2}
\forall s \in [1,2[ \quad\forall\delta>0\quad : \quad\|e\|_{L^{s,m-s-\delta}}<\infty\quad .
\end{equation}
We claim that $v = u \tau$ is H\"older continuous in $B$, provided the bound
\eqref{05} holds with $\varepsilon(m) > 0$ sufficiently small.

Let $B_R(x_0) \subset B$ and let 
\begin{equation}
    P^{-1}dv=df + *dg +h
\end{equation}
be the Hodge decomposition of $P^{-1}dv$ on $B_R(x_0)$, where $f \in H^1_0(B_R(x_0))$
and where $g$ is a co-closed $m-2$-form of class $H^1(B_R(x_0))$ whose restriction to 
the boundary $\p B$ also vanishes, 
and with a harmonic $1$-form $h\in L^2(B_R(x_0))$; see \cite{IwGa} Corollary 10.5.1, p.236,
for the Hodge decomposition of forms in Sobolev Spaces. Similar to \eqref{1.7} we have 
the equations (up to sign, which is of no importance in what follows)
\begin{equation}\label{3.7}
  \begin{split}
    & -\Delta f = - div(P^{-1} \nabla v)= *d\xi \cdot P^{-1}dv - e,\\
    & -\Delta g = \ast d(P^{-1}dv) = \ast (dP^{-1} \wedge dv).
  \end{split}
\end{equation}

Fix a number $1 < p < m/(m-1)$ and let $q > m$ be the conjugate exponent with 
$1/p + 1/q = 1$. Since $f=0$ on $\p B_R(x_0)$, by duality we have 
\begin{equation}\label{3.8}
    ||df||_{L^p} \le C \sup_{\varphi \in W^{1,q}_0(B_R(x_0)); ||\varphi||_{W^{1,q}} \le 1}
    \int_{B_R(x_0)} df \cdot d\varphi \, dx.
\end{equation}
Here and in the following computations all norms refer to the domain $B_R(x_0)$.
Note that $W^{1,q}_0(B_R(x_0)) \hookrightarrow C^{1-m/q}(B_R(x_0))$ and
for all $\varphi \in W^{1,q}_0(B_R(x_0))$ with $||\varphi||_{W^{1,q}} \le 1$
there holds
\begin{equation}\label{3.8a}
  ||\varphi||_{L^{\infty}} \le C R^{1-m/q}||\varphi||_{W^{1,q}} \le C R^{1-m/q}.
\end{equation}

For such $\varphi$ then we can estimate 
\begin{equation}\label{3.9a}
  \begin{split}
    & \int_{B_R(x_0)} df \cdot d\varphi \, dx = - \int_{B_R(x_0)} \Delta f \varphi \, dx\\
    & = \int_{B_R(x_0)} d\xi \wedge P^{-1} dv \varphi - \int_{B_R(x_0)} e \varphi \, dx = I + II
  \end{split}
\end{equation}
as follows. 
Similar to the approach introduced in \cite{CWY}, upon integrating by parts and using 
\cite{CLMS}, Theorem II.1, we have (again up to sign)
\begin{equation}\label{3.9b}
  \begin{split}
  I & = \int_{B_R(x_0)} d\xi \wedge P^{-1} dv \varphi    
    = \int_{B_R(x_0)} d\xi \wedge d(P^{-1} \varphi) (v - \bar{v}_{x_0,R})\\
    & \le C ||d\xi||_{L^2} ||d(P^{-1} \varphi)||_{L^2} [v]_{BMO}\\
    & \le C ||d\xi||^2_{L^2}(||dP||_{L^2} ||\varphi||_{L^{\infty}} + ||d\varphi||_{L^2})  [v]_{BMO}\\
    & \le C R^{m -1-m/q}||d\xi||^2_{L^{2,m-2}}
         (||dP||_{L^{2,m-2}} + ||d\varphi||_{L^q})  [v]_{BMO}\\
    & \le CR^{m/p-1}\varepsilon(m) [v]_{BMO},
  \end{split}
\end{equation}
while \eqref{3.6b}, combined with (\ref{zx2}) and (\ref{3.8a}), for any $\delta>0$ gives the bound
\begin{equation}\label{3.9c}
  \begin{split}
    II & = - \int_{B_R(x_0)} e \varphi \, dx \le ||e||_{L^1(B_R(x_0))} ||\varphi||_{L^{\infty}}\\
       & \le C_\delta R^{m-1-\delta} ||e||_{L^{1,m-1-\delta}}||\varphi||_{L^\infty} \le C_\delta R^{m - m/q-\delta} = C_\delta R^{m/p-\delta}.
  \end{split}
\end{equation}
Hence from (\ref{3.7}) we conclude that for every $\delta>0$ there holds
\begin{equation}\label{3.9}
    ||df||_{L^p} \le CR^{m/p-1}\varepsilon(m) [v]_{BMO} + C_\delta R^{m/p-\delta}.
\end{equation}

Similarly, letting $s$ satisfy $1/2 + 1/q + 1/s = 1$, by H\"older's inequality 
for an arbitrary $m-2$ form $\psi \in W^{1,q}(B_R(x_0),\wedge^{m-2}{\R}^m)$ vanishing 
on $\p B$ and with $||\psi||_{W^{1,q}} \le 1$ we can bound 
\begin{equation}\label{3.10a}
  \begin{split}
    & \int_{B_R(x_0)} dg \cdot d\psi \, dx = - \int_{B_R(x_0)} \Delta g \psi \, dx\\
    & = \int_{B_R(x_0)} dP^{-1} \wedge  dv \psi 
    = \int_{B_R(x_0)} dP^{-1} \wedge d \psi (v - \bar{v}_{x_0,R})\\
    & \le C ||dP||_{L^2} ||d\psi||_{L^q} ||v - \bar{v}_{x_0,R}||_{L^s}
     \le CR^{m/p-1}\varepsilon(m) [v]_{BMO}.
  \end{split}
\end{equation}
By duality, we have
\begin{equation}\label{3.10b}
    ||dg||_{L^p} =\sup_{k\in L^q(B_R(x_0);\wedge^{m-2}{\R}^m);\, \|k\|_q\le 1}
    \int_{B_R(x_0)}dg\cdot k\, .
\end{equation}
Decomposing any $k \in L^q(B_R(x_0);\wedge^{m-2}{\R}^m)$ as $k = d\psi + \ast d\rho + \eta$
with $\psi = 0$ on $\p B$ as in \cite{IwGa}, Corollary 10.5.1, and recalling that
the restriction of $g$ to $\p B_R(x_0)$ vanishes, we then arrive at the estimate
\begin{equation}\label{3.10}
\begin{array}{rl}
    ||dg||_{L^p} &\ds\le C\sup_{\psi \in W^{1,q}(B_R(x_0),\wedge^{m-2}{\R}^m); \, 
                             \|d\psi\|_q\le 1}\int_{B_R(x_0)}dg\cdot  d\psi \\[5mm]
    &\ds\le CR^{m/p-1}\varepsilon(m) [v]_{BMO}.
\end{array}
\end{equation}

From the Campanato estimates for harmonic functions, as in Giaquinta \cite{Gia}, 
proof of Theorem 2.2, p.84 f.,
we thus conclude that for any $r < R$ there holds
\begin{equation}\label{3.11}
  \begin{split}
    \int_{B_r(x_0)} |dv|^p \, dx 
      & \le  \int_{B_r(x_0)} |h|^p \, dx + \int_{B_r(x_0)} (|df|^p+|dg|^p) \, dx\\
      & \le C \left(\frac{r}{R}\right)^m \int_{B_R(x_0)} |h|^p \, dx 
         + \int_{B_r(x_0)} (|df|^p+|dg|^p) \, dx\\
      & \le C  \left(\frac{r}{R}\right)^m \int_{B_R(x_0)} |dv|^p \, dx 
         + \int_{B_R(x_0)} (|df|^p+|dg|^p) \, dx\\
      &  \le C  \left(\frac{r}{R}\right)^m \int_{B_R(x_0)} |dv|^p \, dx 
         + CR^{m-p}\varepsilon(m) [v]^p_{BMO} + C_\delta R^{m-\delta p},
   \end{split}
\end{equation}
for any $\delta>0$.
Set 
\begin{equation*}
  \Phi(x_0,r) = r^{p-m} \int_{B_r(x_0)} |dv|^p \, dx
\end{equation*}
and define
\begin{equation*}
  \Psi(R) =  \sup_{x_0,\, 0 < r < R} \Phi(x_0,r).
\end{equation*}
Then we can bound
\begin{equation*}
  \sup_{x_0}\, [v]^p_{BMO(B_R(x_0))} \le C  \Psi(R),
\end{equation*}
and from \eqref{3.11} we have 
\begin{equation}\label{3.12}
  \begin{split}
     \Phi(x_0,r) & \le C \left(\frac{r}{R}\right)^p \Phi(x_0,R)
         + C\left(\frac{r}{R}\right)^{p-m}\varepsilon(m)  \Psi(R)
         + C_\delta \left(\frac{r}{R}\right)^{p-m} R^{p-\delta p}.
   \end{split}
\end{equation}

Fixing the ratio $r/R = \gamma$ for some number $0 < \gamma <1$ to be 
specified below, we pass to the supremum with respect to 
$x_0 \in B$ and $0 < r = \gamma R < R$.
With a uniform constant $C_1$ from \eqref{3.12} 
for any $R > 0$ we obtain 
\begin{equation}\label{3.13}
  \begin{split}
     \Psi(\gamma R) & \le C_1 \gamma^p(1+ \varepsilon(m)\gamma^{-m}) \Psi(R)
        + C_\delta \gamma^{p-m} R^{p-\delta p},
    \end{split}
\end{equation}
for every $\delta>0$. Now also fix $\delta > 0$ such that $p- p\delta>1$.
Choosing $\gamma$ such that $C_1\gamma^{(p-1)/2} \le 1/2$,
we determine $\varepsilon(m) > 0$ so that $\varepsilon(m) \le \gamma^m$.
With a uniform constant $C$ we then obtain the estimate 
\begin{equation}\label{3.14}
  \begin{split}
     \Psi(\gamma R) & \le \gamma^{(p+1)/2} \Psi(R) + CR^{p-\delta p} 
     \le \gamma^{(p+1)/2} \Psi(R) + CR
   \end{split}
\end{equation}
for all $R \in ]0,1]$.  Finally, for $R = 1$ and any $r \in ]0,\gamma]$, 
letting $k \in \N$ be such that $\gamma^{k+1} < r \le \gamma^k$ and 
iterating as in Giaquinta \cite{Gia}, proof of Lemma 2.1, p.86, we conclude that 
\begin{equation}\label{3.15}
  \begin{split}
    \gamma^{m-p}\Psi(r) & \le \Psi(\gamma^k) \le  \gamma^{k(p+1)/2} \Psi(1) 
     + C\gamma^k(\sum_{j=1}^{\infty}\gamma^{j(p-1)/2}) \le C r.
   \end{split}
\end{equation}
Hence $v \in C^{1/p}(B)$ and therefore also $u \in C^{1/p}(B_{1/2}(0))$, as claimed.

\section{Proof of Lemma~\ref{lm-3.1}}
For the proof of Lemma~\ref{lm-3.1} we follow \cite{MR}, where Uhlenbeck's \cite{Uhl} 
construction of a local Coulomb gauge in Sobolev spaces was generalized to 
Morrey spaces. 
Due to the fact that the space $L^{2,m-2}_1$ defined earlier does not embed into $C^0$, 
the inverse mapping $P\rightarrow P^{-1}$ is not smooth from the space $L^{2,m-2}_1$ into 
itself. 
In order to avoid this difficulty, similar to \cite{Uhl} we first construct the local 
Coulomb gauge under slightly more stringent regularity assumptions.

\begin{lemma}
\label{lm-3.2} There exists $\ep(m,n)>0$ and $C>0$ such that, on $B = B^m$
for every $\al>0$ and every $\Om$ in
 $L^{2,m-2+\al}(B,so(n))$ with 
\begin{equation}\label{4.1}
  ||\Omega||^2_{L^{2,m-2}} \le \varepsilon(n,m)
\end{equation}
there exist $P \in L^{2,m-2+\al}_1(B;SO(n))$ and $\xi \in L^{2,m-2+\al}_1(B,so(n))$ such that 
\begin{equation}\label{3.z00}
  P^{-1}dP +P^{-1}\Omega P = *d\xi \mbox{ on } B, \mbox{ and }\xi=0\mbox{ on }\p B\quad .
\end{equation}
Moreover, $dP$ and $d\xi$ satisfy
\begin{equation}\label{3.z0}
 \begin{split}
   ||dP||^2_{L^{2,m-2+\al}}  + ||d\xi||^2_{L^{2,m-2+\al}} 
   \le C ||\Omega||^2_{L^{2,m-2+\al}} \quad,
  \end{split}
\end{equation}
and
\begin{equation}\label{3.z1}
 \begin{split}
   ||dP||^2_{L^{2,m-2}}  + ||d\xi||^2_{L^{2,m-2}} 
   \le C ||\Omega||^2_{L^{2,m-2}} \le C\varepsilon(n,m).
  \end{split}
\end{equation}
\end{lemma}

\begin{proof}[Proof of Lemma~\ref{lm-3.1}.] 
Let $\Om$ be in $L^{2,m-2}$ and suppose that $\|\Om\|_{L^{2,m-2}}<\ep$ for some number $\ep > 0$
to be fixed below. Although smooth functions are not dense in $L^{2,m-2}$, it is not difficult to
show that the mollified forms $\Om_\delta=\Om\ast\chi_\delta$ obtained from $\Om$ by convoluting
$\Om$ with a standard mollifyer satisfy the uniform estimate 
$\|\Om_\delta\|_{L^{2,m-2}}\le C\|\Om\|_{L^{2,m-2}}$. By choosing
$\ep > 0$ sufficiently small, we can then achieve the uniform bound 
$\|\Om_\delta\|_{L^{2,m-2}}\le\ep(m,n)$, where $\ep(m,n)$ is given in Lemma~\ref{lm-3.2},
to obtain the existence of $\xi_\delta$ and $P_\delta$ satisfying \eqref{3.z00}, \eqref{3.z0} 
and \eqref{3.z1} for $\Om_\delta$ instead of $\Om$. 
The uniform bound given by \eqref{3.z1} permits to pass to the limit
$\delta \to 0$ in \eqref{3.z00}, and the assertion of Lemma~\ref{lm-3.1} follows. 
\end{proof}

\begin{proof}[Proof of Lemma~\ref{lm-3.2}.] 
For $\al>0$ introduce the set
\[
{\mathcal U}^\al_{\ep,C}:=\lf\{
\begin{array}{c}
\Om\in L^{2,m-2+\al}(B^m,so(n));\ \|\Om\|_{L^{2,m-2}}\le \ep, \mbox{ and}\\[3mm]
\mbox{there exist } P\mbox{ and } \xi \mbox{ satisfying }(\ref{3.z00}),\  (\ref{3.z0}),\ (\ref{3.z1})
\end{array}
\rg\}
\]
Since clearly $\Om = 0 \in {\mathcal U}^\al_{\ep,C}$, the set ${\mathcal U}^\al_{\ep,C}$ is not empty.
The proof therefore will be complete once we show that, for $\ep$ small enough and $C$ large enough,
${\mathcal U}^\al_{\ep,C}$ is both open and closed in the star-shaped and hence path-connected set
\[
{\mathcal V}^\al_{\ep}:=\lf\{
\Om\in L^{2,m-2+\al}(B^m,so(n));\ \|\Om\|_{L^{2,m-2}}\le\ep\rg\}.
\]
The proof of closedness is similar to the proof of Lemma~\ref{lm-3.1} given above.
To see that ${\mathcal U}^\al_{\ep,C}$ is open in ${\mathcal V}^\al_{\ep}$, observe that
for $\al>0$ the space $L^{2,m-2+\al}_1$ embeds continuously into $C^0$ and the 
inverse mapping $P\rightarrow P^{-1}$ from the space $L^{2,m-2+\al}(B,SO(n))$ into itself is smooth.
Therefore the argument of \cite{Uhl} can be applied to show that, for sufficiently small $\ep > 0$
and sufficiently large $C$, for every $\Om$ in ${\mathcal U}^\al_{\ep,C}$ there exists 
$\eta_\Om>0$ with the property that for every $\om\in L^{2,m-2+\al}$
satisfying the bound $\|\om\|_{L^{2,m-2+\al}}\le\eta_\Om$ and $\|\Om+\om\|_{L^{2,m-2}}<\ep$ we can find
$\xi_\om$ and $P_\om$ in $L^{2,m-2+\al}_1(B,so(n))$ and
$L^{2,m-2+\al}_1(B,SO(n))$, respectively, satisfying (\ref{3.z00}). The openness
of ${\mathcal U}^\al_{\ep,C}$ may be obtained as in \cite{Riviere} from the following lemma.
This completes the proof. 
\end{proof}
 
\begin{lemma}\label{lm-3.3} 
There exists $\delta>0$ with the following property. 
Suppose that for $\Om \in {\mathcal V}^\al_{\ep}$ there exist 
$\xi \in L^{2,m-2+\al}_1(B^m,so(n))$, $P \in L^{2,m-2+\al}(B^m,SO(n))$ 
satisfying (\ref{3.z00}) and the estimate
\begin{equation}
\label{3.z3}
 \|d\xi\|_{L^{2,m-2}}+\|dP\|_{L^{2,m-2}}\le\delta\quad .
\end{equation}
Then (\ref{3.z0}) and (\ref{3.z1}) hold for some $C$ independent of 
$\Om \in {\mathcal V}^\al_{\ep}$.
\end{lemma}
 
\begin{proof}
In view of (\ref{3.z00}), the $2$-form $\ast\xi$ satisfies
\begin{equation}
\label{3.z4}
\lf\{
\begin{array}{l}
\ds\Delta\ast\xi=dP^{-1}\wedge dP+d(P^{-1}\Om P)\quad\mbox{ in }B ,\\[4mm]
\ds\ast\xi=0\quad\quad\mbox{ on }\p B
\end{array}
\rg.
\end{equation}

We decompose $\ast\xi$ in two forms $\ast\xi=u+v$ solving, respectively,
\begin{equation}
\label{3.z8}
\lf\{
\begin{array}{l}
\ds\Delta u=dP^{-1}\wedge dP\quad\quad\mbox{ in }B\\[4mm]
\ds u=0\quad\quad\mbox{ on }\p B\ ,
\end{array}
\rg.
\end{equation}
and
\begin{equation}
\label{3.z9}
\lf\{
\begin{array}{l}
\ds\Delta v=d(P^{-1}\Om P)\quad\quad\mbox{ in }B\\[4mm]
\ds v=0\quad\quad\mbox{ on }\p B\ .
\end{array}
\rg.
\end{equation}
From \cite{Gia}, Theorem 2.2, 
for $s \in \{2 - \alpha, 2\}$ first we obtain the bound
\begin{equation}
\label{3.z10}
\|dv\|_{L^{2,m-s}(B)}\le C\|\Om\|_{L^{2,m-s}(B)}\ .
\end{equation}
Following the strategy of the proof of Theorem 3.3 in \cite{Gia}, and estimating 
the term on p.48, l.2 in \cite{Gia}
\begin{equation}\label{4.2}
  \begin{split}
   & \|d(u-v)\|^2_{L^2(B_R(x_0)\cap B)} = - \int_{B_R(x_0)\cap B} \Delta (u-v)\wedge (u-v) \\
   & \quad = \int_{B_R(x_0)\cap B} dP^{-1}\wedge dP \wedge (u-v)\\
   & \quad = \int_{B_R(x_0)\cap B} dP^{-1} (P - \bar{P}_{x_0,R}) \wedge d(u-v)\\
   & \quad \le \|d(u-v)\|_{L^2(B_R(x_0)\cap B)}\; \|P\|_{BMO(B)}\; \|dP^{-1}\|_{L^2(B_r(x_0)\cap B)}
  \end{split}
\end{equation}
via \cite{CLMS}, Theorem II.1, likewise for $s \in \{2 - \alpha, 2\}$ we obtain that 
\begin{equation}\label{4.3}
  \begin{split}
   \|du\|_{L^{2,m-s}(B)}& \le C\|P\|_{BMO(B)}\|dP\|_{L^{2,m-s}(B)}\\
   & \le C\|dP\|_{L^{2,m-2}(B)}\|dP\|_{L^{2,m-s}(B)}\le C\delta \|dP\|_{L^{2,m-s}(B)}\ .
  \end{split}
\end{equation}
Combining (\ref{3.z10}) and (\ref{4.3}) we then conclude
\begin{equation}
\label{3.z12}
\|d\xi\|_{L^{2,m-s}(B)}\le C\delta\ \|dP\|_{L^{2,m-s}(B)}+C\ \|\Om\|_{L^{2,m-s}(B)}\quad .
\end{equation}
Moreover, from (\ref{3.z00}) we have
\begin{equation}
\label{3.z13}
\|dP\|_{L^{2,m-s}(B)}\le \|d\xi\|_{L^{2,m-s}(B)}+\|\Om\|_{L^{2,m-s}(B)}\quad .
\end{equation}
Putting the estimates (\ref{3.z12}) and (\ref{3.z13}) together, upon choosing $\delta > 0$ small
enough we then obtain (\ref{3.z0}) and (\ref{3.z1}). The proof is complete.
\end{proof}

\end{document}